\documentclass[12pt,a4paper]{article}
\usepackage{amsmath}
\usepackage{amsthm}
\usepackage{amssymb}
\usepackage{mathrsfs}
\usepackage{enumitem}
\usepackage[left=3cm,right=3cm,top=3cm,bottom=3cm]{geometry}

\usepackage{natbib}
\usepackage[breaklinks=true,colorlinks,ps2pdf,citecolor=blue]{hyperref}
\usepackage[all]{hypcap}
\usepackage{graphicx}

\usepackage{booktabs}
\usepackage{color}

\usepackage[labelfont={bf}]{caption,subfig}
\usepackage{titling}

\thanksmarkseries{arabic}
\date{\today}
\author{Katar\'{\i}na Burclov\'{a}\thanks{katarina.burclova@gmail.com} and Andrej P\'{a}zman\thanks{pazman@fmph.uniba.sk}\\
Comenius University in Bratislava\thanks{Faculty of Mathematics, Physics and Informatics, Mlynsk\'{a} dolina, 842 48 Bratislava, Slovak Republic}
}
\title{Optimum design via $I$-divergence for stable estimation in generalized regression models}

\newtheorem{theorem}{Theorem}

\theoremstyle{definition}
\newtheorem{lemma}{Lemma}

\theoremstyle{definition}
\newtheorem{remark}{Remark}

\begin{document}
\maketitle

\begin{abstract}Optimum designs for parameter estimation in generalized regression models are standardly based on the
Fisher information matrix (cf.~\cite{AFH14} for a recent exposition).
The corresponding optimality criteria are related to the asymptotic
properties of maximal likelihood (ML) estimators in such models. However, in finite sample
experiments there could be problems with identifiability, stability and
uniqueness of the ML estimate, which are not reflected by the information matrices.
In~\cite{PP14} is discussed how to solve some of these estimability issues on the design stage of an
experiment in standard nonlinear regression. Here we want to extend this design methodology to more general
models based on exponential families of distributions (binomial, Poisson,
normal with parametrized variances, etc.). The main tool for that is the
information (or Kullback-Leibler) divergence, which is closely related to the ML estimation. 
\end{abstract} 

\section{Introduction} 
\label{sec:1}

To each design point $x\in \mathcal{X}$, the design space, we associate an
observation $y$ (a~random variable or vector), which is distributed
according to the density of an exponential form 
\begin{equation}
f\left( y\mid \theta ,x\right) =\exp \left\{ -\psi \left( y\right)
+t^\top\left( y\right) \gamma -\kappa \left(\gamma\right) \right\} _{\gamma
=\gamma \left( x,\theta \right) }\;,
\label{eq1}
\end{equation}
with the unknown parameter $\theta $ taking values from a given parameter
space $\Theta \subset \mathbb{R}^p$. This density is taken with respect to a
measure $\nu \left( \cdot\right) $ on $Y$, the sample space of $y$. Usually $Y\subset \mathbb{R}^s$ and $\nu $ is the Lebesgue measure, or $Y$ is finite or
countable and $\nu \left( \left\{ y\right\} \right) =1$ for every $y\in Y$; then $f\left( y\mid \theta ,x\right) $ is simply the probability of $y$.
Well known examples are the one-dimensional normal density 
\[
f\left( y\mid \theta ,x\right) =\frac 1{\left( 2\pi \right) ^{1/2}\sigma
\left( x,\theta \right) }\exp \left\{ -\frac{\left[ y-\mu \left( x,\theta
\right) \right] ^2}{2\sigma ^2\left( x,\theta \right) }\right\} ;\quad y\in 
\mathbb{R} 
\]
or the binomial probability distribution 
\begin{equation}
f\left( y\mid \theta ,x\right) =\binom{n}{y}\pi \left( x,\theta \right) ^y\left( 1-\pi \left( x,\theta \right)
\right) ^{n-y};\quad y\in \left\{ 0,1,\hdots,n\right\}\;. 
\label{eq2}
\end{equation}

Consider an exact design $X=\left( x_1,\hdots,x_N\right)$, where $x_i\in 
\mathcal{X}$ and the observations $y_{x_1},\hdots,y_{x_N}$ are independent. The
ML estimator for $\theta $ is $\hat{\theta}=\arg \max_{\theta
\in \Theta }\sum_{i=1}^N\ln f\left( y_{x_i}\mid \theta ,x_i\right)$. For
large $N$, and under some regularity assumptions, $\hat{\theta}$ is
approximately distributed normally with mean $\theta $ and variance $
M^{-1}\left( X,\theta \right) $, where $M\left( X,\theta \right)
=\sum_{i=1}^N $ $M\left( x_i,\theta \right)$, and $M\left( x,\theta
\right) =E_\theta \left( -\frac{\partial ^2\ln f\left( y\mid \theta
,x\right) }{\partial \theta \partial \theta ^\top}\right) $ is the  elemental
information matrix at $x $ (cf.~\cite{AFH14} for this terminology). Hence within this
asymptotic approximation, a design $X$ is (locally) optimal if it maximizes $
\Phi \left[ M\left( X,\theta^0\right) \right] $, where $\theta^0$ is a
guess for the true (but unknown) value of $\theta $. Here $\Phi \left(\cdot\right) $ stands for $\det^{1/p}(\cdot) $ (or $\ln \det (\cdot)$) in case of $D$-optimality, etc. This is the standard way to
optimize designs in generalized regression models.

Alternatively to the information matrix, we take here for design purposes
the $I$-divergence (the information or Kullback-Leibler divergence, cf.~\cite{K97}), which for any two
points $\theta^0,\theta \in \Theta $ is equal to 
$
I_X\left( \theta^0,\theta \right) =\sum_{i=1}^NI_{x_i}\left( \theta^0,\theta \right) $
with the elemental $I$-divergence defined by 
\begin{equation}
I_x\left( \theta^0,\theta \right) =E_{\theta^0}\left[ \ln \frac{
f\left( y\mid \theta^0,x\right) }{f\left( y\mid \theta ,x\right) }\right]\;.
\label{eq3}
\end{equation}

As is well known, $I_x\left( \theta^0,\theta \right) \geq 0$, and it is
equal to zero if and only if $f\left( y\mid \theta^0,x\right) =f\left(
y\mid \theta ,x\right)$. In general, the
$I$-divergence measures well the sensitivity of the data $y$ to the shift
of the parameter from the value $\theta^0$ to the value $\theta $, even
when $\theta $ and $\theta^0$ are distant, while the information matrix $M\left( x,\theta^0\right) $ is doing essentially the same, but only for $\theta $ which is close to $\theta^0$ (see Sect.~\ref{sec:3}). Hence the
$I$-divergence may allow a better characterization of the statistical
properties of the model than the information matrix. An important fact is
also that we can compute it easily (avoiding integrals) in models given by~(\ref{eq1}). Notice that for the normal model with $\sigma ^2\left( x,\theta \right)\equiv 1$, one has $2I_X\left( \theta^0,\theta \right)
=\sum_{i=1}^N\left[ \mu \left( x_i,\theta^0\right) -\mu \left( x_i,\theta
\right) \right] ^2$, an expression, which is largely used in~\cite{PP14}. 

We note that in~\cite{LTT07} the $I$-divergence has been used for design purposes for model discrimination, which, however, is a different aim.

\section{Basic Properties of Model~(\ref{eq1})} 
\label{sec:2}
It is clear that $t\left( y\right) $ is a sufficient statistics in model~(\ref{eq1}), so we can suppose, at least in theory, that we observe $t\left(y\right) $ instead of $y$. Denote by $\eta \left( x,\theta \right) $ its
mean. For  $\gamma =\gamma ( x,\theta ) $ we have
\begin{equation}
\eta \left( x,\theta \right) =\int_Y t( y) \exp \left\{ -\psi
( y) +t^\top( y) \gamma -\kappa \left(\gamma \right)\right\}_{\gamma =\gamma \left( x,\theta \right)} d\nu ( y)=\left[ \frac{\partial \kappa \left( \gamma \right) }{\partial \gamma }\right] _{\gamma =\gamma \left( x,\theta \right) }\;.
\label{eq4}
\end{equation}
To be able to do this derivative at any $\gamma$, we suppose that the set
$\bigl\{ \gamma :\int_Y\exp \bigl\{ -\psi \left( y\right)+t^\top\left( y\right)\gamma \bigr\} d\nu \left( y\right) <\infty \bigr\} $
is open. Then the model~(\ref{eq1}) is called regular, and regular models are
standard in applications. Taking the second order derivative in~(\ref{eq4}) we obtain $\frac{\partial ^2\kappa \left( \gamma \right) }{\partial \gamma \partial
\gamma ^\top}=Var_\gamma \left[ t\left( y\right) \right]$, which for $\gamma
=\gamma \left( x,\theta \right) $ will be denoted by $\Sigma \left( x,\theta
\right) $. By a reduction of the linearly dependent components of the vector 
$t\left( y\right) $ one can always achieve that $\Sigma \left( x,\theta
\right) $ is nonsingular, and we obtain from~(\ref{eq4}) that $\frac{\partial \eta \left( x,\theta \right) }{\partial \theta ^\top}=\Sigma\left( x,\theta \right) \frac{\partial \gamma \left( x,\theta \right) }{\partial \theta ^\top}$.

The functions $\theta \in \Theta \rightarrow \eta \left( x,\theta \right) $
(the mean-value function) and $\theta \in \Theta \rightarrow \gamma \left(
x,\theta \right) $ (the canonical function) are useful dual representations
of the family of densities~(\ref{eq1}) (cf.~\cite{E78}). 
From~(\ref{eq1}) and~(\ref{eq4}) it follows that $E_{x,\theta }\left[ \frac{\partial \ln f\left( y\mid x,\theta \right) }{\partial \theta }\right] =0$, and
consequently the elemental information matrix is equal to 
\begin{equation}
M\left( x,\theta \right)=Var_{x,\theta }\left[ \frac{\partial \ln f\left(
y\mid x,\theta \right) }{\partial \theta }\right] 
=\frac{\partial \eta^\top\left( x,\theta \right) }{\partial \theta }\Sigma^{-1} \left( x,\theta \right) 
\frac{\partial \eta \left( x,\theta \right) }{\partial \theta ^\top} \;.
\label{infomatica}
\end{equation}
The elemental $I$-divergence is, according to~(\ref{eq3}) and~(\ref{eq1}),
\[
I_x\left( \theta^0,\theta \right) =\eta^\top\left( x,\theta^0\right)
\left[ \gamma \left( x,\theta^0\right) -\gamma \left( x,\theta \right)
\right] +\kappa \left( \gamma \left( x,\theta \right) \right) -\kappa \left(
\gamma \left( x,\theta^0\right) \right)\;.\]
For more details on exponential families see~\cite{B86}.

\section{Variability, Stability and $I$-Divergence} 
\label{sec:3}
In this section we consider observations according to an ``exact" design $X=\left(x_1,\hdots,x_N\right)$.
 The joint density of $\mathbf{y}=\left( y_{x_1},\hdots,y_{x_N}\right) ^\top$ is equal to $\widetilde{f}\left( \mathbf{y\mid }\theta \right) =\prod_{i =1}^Nf\left( y_{x_i}\mid \theta ,x_i\right) $.

The variability of the ML estimate $\hat{\theta}$ in the neighborhood of $\bar{\theta}$, the true value of $\theta $, is well expressed by the information matrix 
$M\left( X,\bar{\theta}\right) $, since its inverse is the asymptotic
variance matrix of $\hat{\theta}$. But the same can be achieved by the
$I$-divergence, since we have in model (\ref{eq1}) the following property of the $I$-divergence.

\begin{lemma} If for any $x\in \mathcal{X}$ the
third order derivatives of $I_x\left( \bar{\theta},\theta \right) $ with respect to $\theta $ are bounded on a neighborhood
of $\bar{\theta}$, then we have  
\begin{equation}
I_X\left( \bar{\theta},\theta \right) =\frac{1}{2}\left( \theta -\bar{\theta}\right) ^\top M\left( X,\bar{\theta}\right) \left( \theta -\bar{\theta}\right)
+o\left( \left\| \theta -\bar{\theta}\right\| ^2\right)\;.
\label{eq5}
\end{equation}
\label{prop1}
\end{lemma}

It is sufficient to prove this equality for the elemental $I$-divergence and
elemental information matrix. We have $I_x\left( \bar{\theta},\bar{\theta}\right) =0,\frac{\partial I_x\left( \bar{\theta},\theta \right) }{\partial
\theta }\mid _{\theta =\bar{\theta}}=0,\frac{\partial ^2I_x\left( \bar{\theta%
},\theta \right) }{\partial \theta \partial \theta ^\top}\mid _{\theta =\bar{\theta}}=M\left( x,\bar{\theta}\right)$, so by the Tylor formula we obtain~(\ref{eq5}). 

On the other hand, we can have important instabilities of the ML estimate $\hat{\theta}$ when, with a
large probability, $\ln \widetilde{f}\left( \mathbf{y\mid }\theta \right) $ is close
to $\ln \widetilde{f}\left( \mathbf{y\mid }\bar{\theta}\right) $ for a point $\theta$ distant from $\bar{\theta}$. However, at the design stage we do not know 
the value of $\mathbf{y}$, so we cannot predict the value of the difference $\ln \widetilde{f}\left( \mathbf{y\mid }\theta \right) -\ln \widetilde{f}\left( 
\mathbf{y\mid }\bar{\theta}\right)$. But we can predict at least its mean.

\begin{lemma}For any $\theta \in \Theta $ we have 
$E_{\bar{\theta}}\left\{ \ln \widetilde{f}\left( \mathbf{y\mid }\bar{\theta}\right) -\ln \widetilde{f}\left( \mathbf{y\mid }\theta \right) \right\}
=I_X\left( \bar{\theta},\theta \right)$.
\label{prop2}
\end{lemma}
This equality is evident from~(\ref{eq3}).

As a consequence of Lemmas~\ref{prop1} and~\ref{prop2} we have that the $I$-divergence $I_X\left( \bar{\theta},\theta \right) $ can express simultaneously both: the  variability of the ML estimate $\hat{\theta}$ in a neighborhood of $\bar{\theta}$ and the danger of instability of the ML estimation due to the
possibility of  ``false" estimates which are very distant from the true
value $\bar{\theta}$.

\section{Extended Optimality Criteria} 
\label{sec:4} 
According to the principle mentioned in Sect.~\ref{sec:3}, the design based on the
$I$-divergence should minimize the variability of $\hat{\theta}$ (related to
the information matrix) and protect against instabilities coming from a
value $\theta $ which is distant from the true value $\bar{\theta}$. This
requirement can be well reflected by extended optimality criteria (cf.~\cite{PP14}
for the classical nonlinear regression). To any design (design measure or
approximate design) $\xi $ on a finite design space $\mathcal{X}$ we define
the extended criteria in a form 
\begin{equation}
\phi _{ext}\left( \xi ,\theta^0\right) =\min_{\theta \in \Theta}\sum_{x\in \mathcal{X}}\left\{ 2I_x\left( \theta^0,\theta \right) \left[\frac 1{\rho^2\left( \theta^0,\theta \right) }+K\right] \right\} \xi
\left( x\right)\;,
\label{eq6}
\end{equation}
where $K\geq 0$ is a tuning constant chosen in advance, $\theta^0$ is a guess
for the (unknown) value of ${\theta}$, $\rho \left( \theta^0,\theta
\right) $ is a distance measure (a norm or a pseudonorm) not depending on
the design $\xi $. When $\rho \left( \theta^0,\theta \right) =\left\|
\theta^0-\theta \right\|$, the Euclidean norm, we have the extended
$E$-optimality criterion, denoted by $\phi _{eE}\left( \xi ,\theta^0\right)$. When $\rho \left( \theta^0,\theta \right) =\bigl| h\left( \theta^0\right) -h\left( \theta \right) \bigr|$, with $h\left( \theta \right)\in\mathbb{R} $,
a given function of $\theta$, we have the extended $c$-optimality criterion,
denoted by $\phi _{ec}\left( \xi ,\theta^0\right)$. When $\rho \left(
\theta^0,\theta \right) =\max_{x\in \mathcal{X}}\left| \alpha \left(
x,\theta^0\right) -\alpha \left( x,\theta \right) \right|$, with $x\in 
\mathcal{X\rightarrow }\alpha \left( x,\theta \right) \in \mathbb{R}$, the regression
function of interest, we have the extended $G$-optimality criterion, denoted $\phi_{eG}\left( \xi ,\theta^0\right)$. Notice that usually $\alpha \left(
x,\theta \right) $ is equal to $\gamma \left( x,\theta \right) $ or to $\eta
\left( x,\theta \right) $, but not necessarily. The names of the criteria
are justified by the following theorem.

\begin{theorem}
Let $\mathcal{B}\left( \theta^0,\delta \right) $ be a
ball centred at $\theta^0$ with the diameter $\delta $, and $M\left(\theta^0,\xi\right)=\sum_{x\in\mathcal{X}}M\left(x,\theta\right)\xi(x)$. Then 
\[\lim_{\delta \rightarrow 0}\min_{\theta \in \mathcal{B}\left( \theta^0,\delta \right) }\sum_{x\in \mathcal{X}}2I_x\left( \theta^0,\theta
\right) \left[ \frac 1{\rho ^2\left( \theta^0,\theta \right) }+K\right]
\xi \left( x\right) \]
is equal to
\begin{itemize}
\item $\lambda _{\min }\left[ M\left( \theta^0,\xi \right) \right] $, the minimal eigenvalue, in case that $\rho \left( \theta^0,\theta \right) =\left\| \theta^0-\theta \right\|$,
\item $\left[ c^\top M^{-}\left( \theta^0,\xi \right) c\right] ^{-1}$
 if $c=\left[ \frac{\partial h\left( \theta \right) }{\partial
\theta }\right] _{\theta^0}$ is in the range of $M\left( \theta^0,\xi \right)$, or zero if it is not in this range, in case of 
$\rho \left( \theta^0,\theta \right) =\left| h\left( \theta^0\right)
-h\left( \theta \right) \right| $,
\item $\left[ \max_{x\in \mathcal{X}}f^\top\left( x\right) M^{-1}\left(
\theta^0,\xi \right) f\left( x\right) \right] ^{-1}$ with $f\left( x\right) =\left[ \frac{\partial \alpha \left( x,\theta \right) }{\partial \theta }\right] _{\theta^0}$ if $M\left( \theta
^0,\xi \right) $ is nonsingular and if  $\rho \left(
\theta^0,\theta \right) =\max_{x\in \mathcal{X}}\left| \alpha \left(x,\theta^0\right) -\alpha \left( x,\theta \right) \right| $.
\end{itemize}
When the model is normal, linear, with unit variances of
observations, $h\left(\theta\right)$ is linear function of $\theta$ and $\alpha\left(x,\theta\right)=\eta\left(x,\theta\right)$ is linear, then $\phi _{ext}\left( \xi ,\theta^0\right) $
coincides with the corresponding well known optimality criterion in linear
models.
\label{prop3}
\end{theorem}
\begin{proof}The proof follows from Lemma~\ref{prop1}  and from known
expressions: $\lambda _{\min }\left[ M\right] =\min_u\frac{u^\top Mu}{u^\top u}$,
 $c^\top M^{-}c=\max_{u:\; Mu\neq 0}\frac{\left( c^\top u\right) ^2}{u^\top Mu}$ (if $c$ is in the range of $M$).
In a normal linear model with unit variances we
use that $\gamma \left( x,\theta \right)=\eta \left( x,\theta \right)=f^\top \left( x\right) \theta $ and $\sum_{x\in \mathcal{X}}2I_x\left( \theta^0,\theta \right) \xi\left( x\right) =\left( \theta -\theta^0\right) ^\top\left[ \sum_{x\in \mathcal{X}}f\left( x\right) f^\top\left(
x\right) \xi \left( x\right) \right] \left( \theta -\theta^0\right)$.
\end{proof}

If $K$ is chosen very large, then the extended criterion gives the same optimum design as the classical criterion. On the other hand, when $K$ is very small, then $1/\rho ^2\left( \theta ^0,\theta \right) $ is the dominating term in~(\ref{eq6}), and
we reject designs even  with non-important instabilities at points $\theta $ distant
from $\theta ^0$.

\begin{figure}[!ht]
     \subfloat[\label{fig:gama}]{%
 \includegraphics[width=0.3\textwidth]{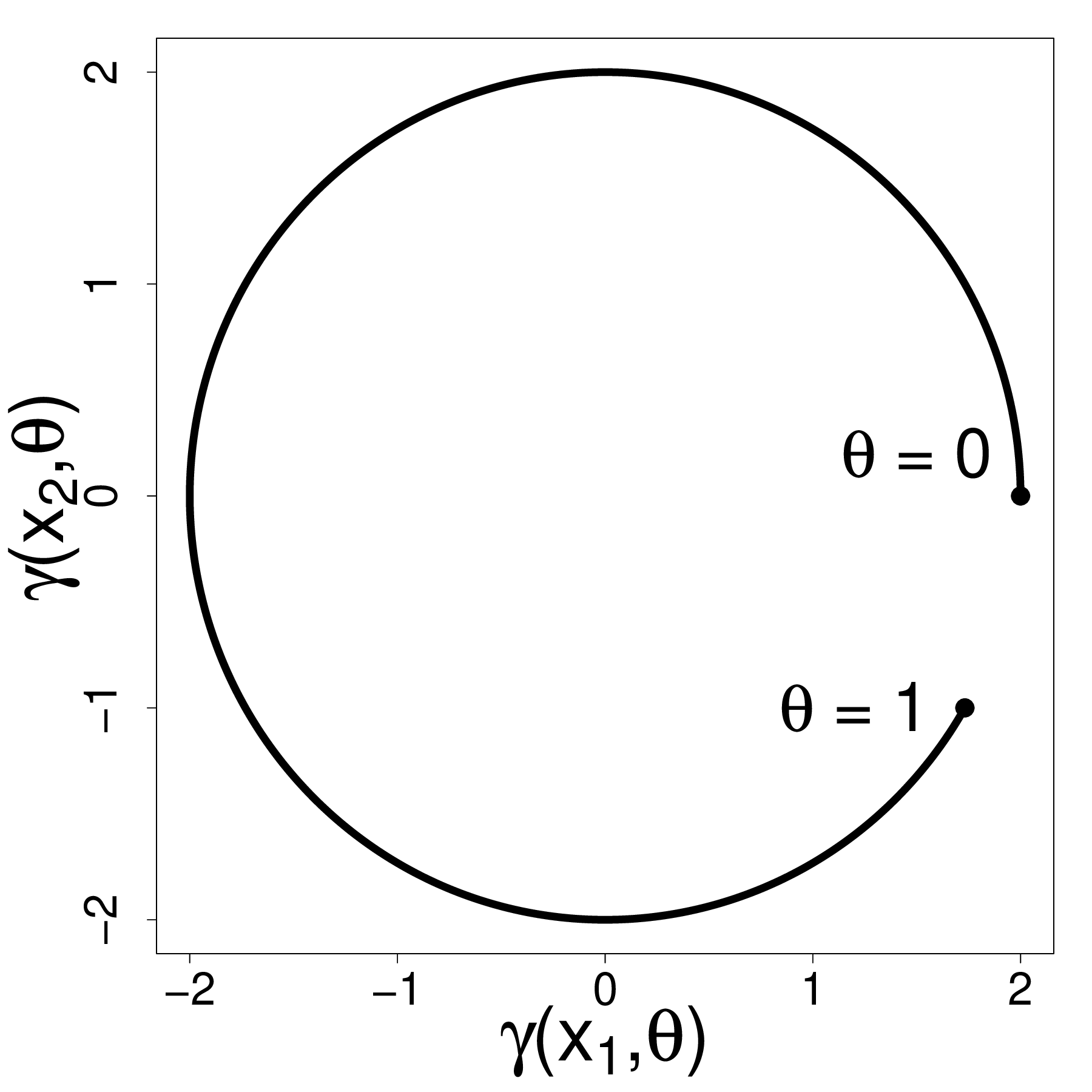}
    }
    \hfill
    \subfloat[\label{fig:eta}]{%
     \includegraphics[width=0.3\textwidth]{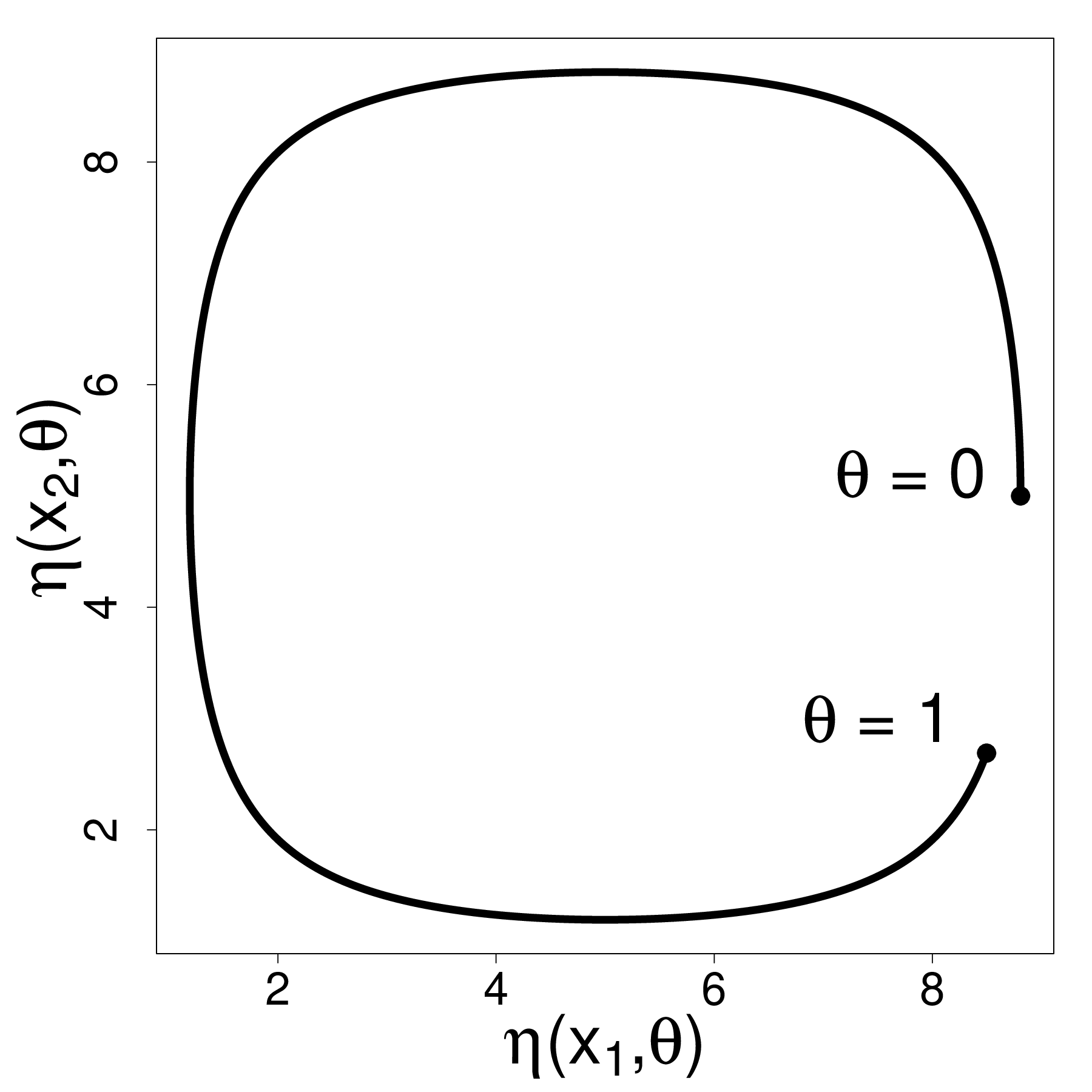}
    }
    \hfill
    \subfloat[ \label{fig:H}]{%
   \includegraphics[width=0.35\textwidth]{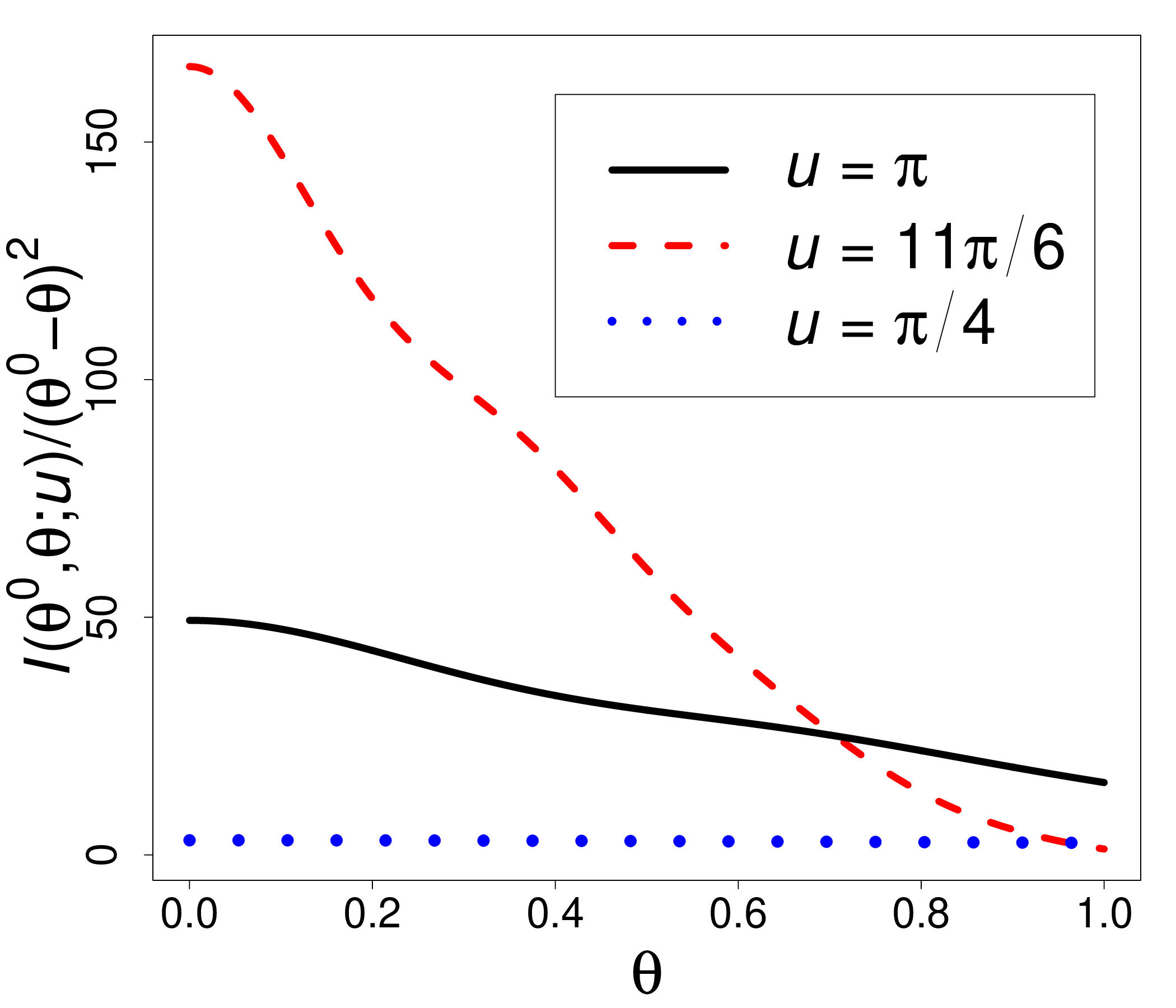}
   }
   \caption{\textbf{(a)} the canonical surface, \textbf{(b)} the expectation surface, \textbf{(c)} $I\left( \theta ^0,\theta ;u\right) /\left( \theta -\theta
^0\right) ^2$}
   \label{fig:dummy}
  \end{figure}

\begin{remark} We can write~(\ref{eq6}) in a form 
\begin{equation}
\phi _{ext}\left( \xi ,\theta^0\right) =\min_{\theta \in \Theta}\sum_{x\in \mathcal{X}}H\left( x,\theta^0,\theta \right) \xi \left(x\right)\;,
\label{eq7}
\end{equation}
with an adequately chosen $H\left( x,\theta^0,\theta \right) $.
It follows that the function $\xi \rightarrow \phi _{ext}\left( \xi,\theta^0\right) $ is concave, hence it has a directional derivative, and the
``equivalence theorem" can be formulated, exactly as in~\cite{PP14}.
\end{remark}
\begin{remark} In the case that the design space $\mathcal{X}$ is finite, $\mathcal{X}=\left\{x^1,\hdots,x^k\right\} $, we can consider the task of computing the optimum design, $\xi^* =\arg \max_\xi \phi _{ext}\left(\xi ,\theta^0\right) $ as an ``infinite-dimensional" linear programming (LP) problem. Namely, we have to
find the vector $\left(t,\xi \left( x^1\right) ,\hdots,\xi \left( x^k\right) \right)$, which maximizes $t$ under the linear constraints 
\begin{eqnarray*}
\sum_{x\in \mathcal{X}}H\left( x,\theta^0,\theta \right) \xi \left(
x\right)  &\geq &t\;\text{for every }\theta \in \Theta\; , \\
\sum_{x\in \mathcal{X}}\xi \left( x\right)  &=&1,\;\xi \left( x\right)  \geq 0\;\text{for every }x\in \mathcal{X}\;.
\end{eqnarray*}
Like in~\cite{PP14}, this can be approximated by iterative LP problems, including a stopping rule, which may be even  more practical than the classical
``equivalence theorem" (see the Numerical example below).
\label{rem2}
\end{remark}
\paragraph{Illustrative Example} 
Consider the binomial model~(\ref{eq2}), which can be written in the exponential
form~(\ref{eq1}): 
$$
f\left( y\mid \theta ,x\right) =\exp \left\{ \ln \binom{n}{y} +y\gamma \left( x,\theta \right) -n\ln \left( 1+e^{\gamma \left(x,\theta \right) }\right)\right\} 
$$
with $\gamma \left( x,\theta \right) =\ln \left[ \pi \left( x,\theta
\right) /\left( 1-\pi \left( x,\theta \right) \right) \right]$, and with
the mean of $y=t(y)$ equal to $\eta \left( x,\theta \right) =n\pi \left(
x,\theta \right) =ne^{\gamma \left( x,\theta \right) }/\left( 1+e^{\gamma
\left( x,\theta \right) }\right) $ (the logistic function). In the example
we took $n=10$, and we considered the regression model (similar to that in~\cite{PP14}) 
$$
\gamma \left( x,\theta \right) =2\cos \left( t-u\theta \right)
;\;x=\left( t,u\right) ^\top\;,
$$
with two observations, one at $x_1=\left( 0,u\right) ^\top$ and the second at $x_2=\left( \pi /2,u\right) ^\top$, where $u\in \left[0,\frac{11}{6}\pi\right]$ is to be
chosen optimally for the estimation of the unknown parameter $\theta \in
[0,1]$. For the case of $u=\frac{11}{6}\pi $ we can see the circular
``canonical surface" $\left\{ \left( \gamma \left( x_1,\theta \right),\gamma \left( x_2,\theta \right) \right) ^\top;\;\theta \in [0,1]\right\} $ in Fig.~\ref{fig:gama}
and the ``expectation surface" $\left\{ \left( \eta \left( x_1,\theta
\right) ,\eta \left( x_2,\theta \right) \right) ^\top;\;\theta \in
[0,1]\right\} $ in Fig.~\ref{fig:eta} (which is no more circular due to the nonlinearity of the
logistic function). The information ``matrix" $M_u\left(\theta \right)\equiv M\left(x_1,\theta\right)+M\left(x_2,\theta\right) $
is computed according to~(\ref{infomatica}), and for $\theta ^0=0$ it is equal to $M_u\left(\theta ^0\right) =nu^2$. It follows that the classical locally optimal design maximizing $M_u\left(\theta ^0\right) $ is obtained when $u=\frac{11}{6}\pi $. We see even from Fig.~\ref{fig:gama} and Fig.~\ref{fig:eta} that under this design, the ML estimate $\hat{\theta}
\left( y\right) $ can be, with a large
probability, in the neighborhood of $\theta =1$, hence the estimator is instable when $\theta ^0=\bar{\theta}=0$.

On the other hand, take the $I$-divergence $I\left( \theta ^0,\theta ;u\right)
\equiv I_{x_1}\left( \theta ^0,\theta \right) +I_{x_2}\left( \theta ^0,\theta
\right) $ with  
\begin{equation}
I_x\left( \theta ^0,\theta \right) =n\left[ \pi \left( x,\theta ^0\right)
\ln \frac{\pi \left( x,\theta ^0\right) }{\pi \left( x,\theta \right) }
+\left( 1-\pi \left( x,\theta ^0\right) \right) \ln \frac{1-\pi \left(
x,\theta ^0\right) }{1-\pi \left( x,\theta \right) }\right] 
\label{eq8}
\end{equation}
and consider the extended criterion $\phi _{ext}\left( u,\theta ^0\right)
=\min_{\theta \in [0,1]}I\left( \theta ^0,\theta ;u\right) /\left( \theta
-\theta ^0\right) ^2$. Numerical computation gives that $\arg \max_{u\in \left[0,\frac{11}{6}\pi\right]}\phi _{ext}\left( u,\theta ^0\right)\doteq\pi $, and for this
choice of $u$ the probability of a false $\hat{\theta}\left( y\right) $ is
negligible, because then $\left(\gamma\left(x_1,\theta\right),\gamma\left(x_2,\theta\right)\right)^\top$ are for $\theta=0$ and $\theta=1$ as far as possible, the same holds for $\eta\left(x,\theta\right)$. We took here the tuning constant $K=0$, otherwise, the optimal $u$ would be between $\pi$ and $\frac{11}{6}\pi$. In Fig.~\ref{fig:H} we present the dependence of $I\left( \theta ^0,\theta ;u\right) /\left( \theta -\theta
^0\right) ^2$ on $\theta $ for different values of $u$.
\paragraph{Numerical Example} The aim of this example is to show that in the case that the design space $\mathcal{X}=\left\lbrace x^1,\hdots,x^k\right\rbrace$ is finite, we can compute the extended $E$-optimum design by using the LP---see Remark~\ref{rem2}. 
The systematic part of the considered model is taken similar as in~\cite{PP14},
Example~2, i.e. the mean of $y_x$, observed at the design point 
$x=\left(x_1,x_2\right)^\top$ is equal to
\begin{equation}
\eta \left( x,\theta \right) =n\pi\left(x,\theta\right)=\frac{n}{6}\left(1+\theta_1x_1+\theta_1^3\left(1-x_1\right)+\theta_2x_2+\theta_2^2\left(1-x_2\right)\right)\;.
\label{eq9}
\end{equation}
However, the error structure is quite different. Instead of an error
component not depending on $\theta =\left(\theta_1,\theta_2\right)^\top$, we have now a binomial model with $y_x$
distributed according to~(\ref{eq2}), with $n=10$ and $\pi\left(x,\theta\right)$ from~(\ref{eq9}). Consequently, in the extended
criterion~(\ref{eq6}) we use the binomial $I$-divergence~(\ref{eq8}) and $\rho ^2\left( \theta^0,\theta \right) =\left\| \theta ^0-\theta \right\| ^2$. We choose $\theta^0=\left(1/8,1/8\right)^\top$, $x\in\mathcal{X}=\{0,0.1,\hdots,0.9,1\}^2$, $\theta \in \Theta=[-1,1]\times[0,2]$. 

The below used iterative algorithm follows the lines of~\cite{PP14}.
\begin{enumerate}[start=0]
\item {Take any vector $\left(\xi^{(0)}\left(x^1\right),\hdots,\xi^{(0)}\left(x^k\right)\right)$ such that $\sum_{x\in\mathcal{X}}\xi^{(0)}(x)=1$ and $\xi^{(0)}(x)\geq 0\;\forall\;x\in\mathcal{X}$, choose $\epsilon>0$, set $\Theta^{(0)}=\emptyset $ and $n=0$. Construct a finite grid $\mathcal{G}^{(0)}$ in $\Theta$.
\item Set $\Theta^{(n+1)}= \Theta^{(n)} \cup \left\lbrace\theta^{(n+1)}\right\rbrace$, where $\theta^{(n+1)}=\arg\min_{\theta\in\Theta}\sum_{x\in\mathcal{X}}H\left(x,\theta^0,\theta\right)\xi^{(n)}(x)$ is computed as follows:}
\begin{enumerate}
\item {Compute $\hat{\theta}^{(n+1)}=\arg\min_{\theta\in\mathcal{G}^{(n)}}\sum_{x\in\mathcal{X}}H\left(x,\theta^0,\theta\right)\xi^{(n)}(x)$.}
\item {Perform local minimization over $\Theta$ initialized at $\hat{\theta}^{(n+1)}$, denote by $\theta^{(n+1)}$ the solution and set $\mathcal{G}^{(n+1)}=\mathcal{G}^{(n)}\cup\left\{\theta^{(n+1)}\right\}$.}
\end{enumerate}
\item {Use the LP solver to find $\left(t^{(n+1)},\xi^{(n+1)}\left(x^1\right),\hdots,\xi^{(n+1)}\left(x^k\right)\right)$ so to maximize $t^{(n+1)}$ satisfying the constraints:}
 \begin{itemize}
 \item {$t^{(n+1)}>0,\; \xi^{(n+1)}(x)\geq 0 \; \forall \; x\in\mathcal{X} ,\; \sum_{x\in\mathcal{X}}\xi^{(n+1)}(x)=1$, }
 \item {$\sum_{x\in\mathcal{X}} H\left(x,\theta^0,\theta\right)\xi^{(n+1)}(x)\geq t^{(n+1)}\; \forall \theta\in\Theta^{(n+1)}$.}
 \end{itemize}
 \item {Set $\Delta^{(n+1)}=t^{(n+1)}-\phi_{ext}\left(\xi^{(n+1)},\theta^0\right)$, if $\Delta^{(n+1)}<\epsilon$, take $\xi^{(n+1)}$ as an $\epsilon$-optimal design and stop, or else $n\leftarrow n+1$ and continue by step 1.}
\end{enumerate}
The computations were performed in Matlab on a bi-processor PC (3.10 Ghz)
equipped with 6GB of RAM and with 64 bits Windows 8.1. LP problems were solved
with interior point method.
When the grid $\mathcal{G}^{(0)}$ was taken as a random latin hypercube design with 10000 points renormalized to $\Theta$, $\epsilon=10^{-10}$, and $\xi^{(0)}$ put mass uniformly to each $x$ in $\mathcal{X}$, the algorithm stopped for $K=0$ after 14 iterations and for $K=10^6$ after 20 iterations requiring 15 and 17\,s. The obtained numerical results are summarized in Table~\ref{Tab}.
\begin{table}[!ht]
 \caption{The $\phi_{eE}$ optimal designs for $K=0$ and $K=10^6$; the values of  $\phi_{eE}$ at $\xi^*$ according to~(\ref{eq6}) for $K=0$ and $K=10^6$}
 \label{Tab}     
\begin{tabular}{|c|c|c|c|}
\hline\noalign{\smallskip}
 $K$ & $\xi^*$ & $\phi_{eE}\left(\xi^*,\theta^0\right)$ for $ K=0$  & $\phi_{eE}\left(\xi^*,\theta^0\right)$ for $ K=10^6$ \\
\noalign{\smallskip}\hline\noalign{\smallskip}
$0$& $\begin{Bmatrix}
    (0,0)^\top&
    (0,1)^\top&
   (1,1)^\top\\
    0.3464&0.0281&0.6255
    \end{Bmatrix}$&0.0215&0.0365\\
 $10^6$& $\begin{Bmatrix}
    (1,0)^\top&
     (0,1)^\top\\
    0.4921&0.5079
    \end{Bmatrix}$&$2.17\times 10^{-9}$&$0.6666$\\
\noalign{\smallskip}\hline\noalign{\smallskip}
\end{tabular}
\end{table}

For $K=10^6$ we obtained almost $E$-optimal design. On the other hand, under this ``$E$-optimal" design, the value of the criterion $\phi_{eE}$ with $K=0$  is small, which indicates instabilities in the model.  

\paragraph{Acknowledgement}
The authors thank Slovak Grant Agency VEGA, Grant No. 1/0163/13, for financial support.

\bibliographystyle{apalike}
\bibliography{references}

\begin{thebibliography}{}

\bibitem[Atkinson et~al., 2014]{AFH14}
Atkinson, A.~C., Fedorov, V.~V., Herzberg, A.~M., and Zhang, R. (2014).
\newblock Elemental information matrices and optimal experimental design for
  generalized regression models.
\newblock {\em J. Stat. Plan. Inference}, 144:81--91.

\bibitem[Brown, 1986]{B86}
Brown, L.~D. (1986).
\newblock {\em Fundamentals of Statistical Exponential Families with
  Applications in Statistical Decision Theory}.
\newblock IMS Lecture Notes Monogr. Ser., Volume 9. Institute of Mathematical
  Statistics, Hayward, CA.

\bibitem[Efron, 1978]{E78}
Efron, B. (1978).
\newblock The geometry of exponential families.
\newblock {\em Ann. Stat.}, 6(2):362--376.

\bibitem[Kullback, 1997]{K97}
Kullback, S. (1997).
\newblock {\em Information Theory and Statistics}.
\newblock Dover Publications, Inc., Mineola, N.Y.

\bibitem[L\'opez-Fidalgo et~al., 2007]{LTT07}
L\'opez-Fidalgo, J., Tommasi, C., and Trandafir, P.~C. (2007).
\newblock An optimal experimental design criterion for discriminating between
  non-normal models.
\newblock {\em J. R. Stat. Soc. B}, 69(2):231--242.

\bibitem[P\'azman and Pronzato, 2014]{PP14}
P\'azman, A. and Pronzato, L. (2014).
\newblock Optimum design accounting for the global nonlinear behavior of the
  model.
\newblock {\em Ann. Stat.}, 42(4):1426--1451.

\end{thebibliography}
\end{document}